\documentclass[12pt,reqno]{amsart}
\usepackage{amssymb}
\usepackage{latexsym}
\usepackage{eucal}

\theoremstyle{plain}
\newtheorem{lemma}{Lemma}[section]
\newtheorem{theorem}[lemma]{Theorem}
\newtheorem{proposition}[lemma]{Proposition}
\newtheorem{corollary}[lemma]{Corollary}
\newtheorem{claim}{Claim}
\newtheorem*{sclaim}{Claim}

\newtheorem*{tha}{Theorem~\ref{T:FinPresEquiv}}
\newtheorem*{thb}{Theorem~\ref{T:SubmFP}}
\newtheorem*{thc}{Theorem~\ref{T:CohSubQ}}
\newtheorem*{thd}{Theorem~\ref{T:FinPres}}

\theoremstyle{definition}
\newtheorem{definition}[lemma]{Definition}

\newtheorem{example}[lemma]{Example}
\newtheorem{problem}{Problem}
\newtheorem*{remark}{Remark}

\numberwithin{equation}{section}

\newcommand{\qedc}{{\qed}~{\rm Claim~{\theclaim}.}}
\newcommand{\sqedc}{{\qed}~{\rm Claim.}} 

\newenvironment{cproof} {\begin{proof}[Proof of Claim.]}
{\qedc\renewcommand{\qed}{}\end{proof}}

\newenvironment{scproof} {\begin{proof}[Proof of Claim.]}
{\sqedc\renewcommand{\qed}{}\end{proof}}

\DeclareMathOperator{\FA}{FA}
\DeclareMathOperator{\PL}{PL}
\DeclareMathOperator{\Conv}{Conv}

\newcommand{\res}{\mathbin{\restriction}}

\newcommand{\poag}{partially ordered abelian group}
\newcommand{\po}{\emph{po}}
\newcommand{\povs}{partially ordered vector space}
\newcommand{\por}{\emph{po}-ring}
\newcommand{\alo}{abelian lattice-ordered}
\newcommand{\cpc}{convex polyhedral cone}
\newcommand{\porm}[1]{partially ordered right $#1$-module}

\newcommand{\NN}{\mathbb{N}}
\newcommand{\ZZ}{\mathbb{Z}}
\newcommand{\QQ}{\mathbb{Q}}
\newcommand{\RR}{\mathbb{R}}

\newcommand{\bA}{\boldsymbol{A}}
\newcommand{\bB}{\boldsymbol{B}}
\newcommand{\bC}{\boldsymbol{C}}
\newcommand{\bD}{\boldsymbol{D}}
\newcommand{\bF}{\boldsymbol{F}}

\newcommand{\bR}{\boldsymbol{R}}
\newcommand{\bS}{\boldsymbol{S}}
\newcommand{\bU}{\boldsymbol{U}}
\newcommand{\bG}{\boldsymbol{G}}
\newcommand{\bH}{\boldsymbol{H}}
\newcommand{\bE}{\boldsymbol{E}}
\newcommand{\bK}{\boldsymbol{K}}

\renewcommand{\SS}{\mathcal{S}}
\newcommand{\TT}{\mathcal{T}}

\newcommand{\MM}{\mathfrak{M}}

\def\vv<#1>{\langle#1\rangle}
\def\fs(#1){\{1,\ldots,#1\}}

\begin{document}

\title[Coherent ordered rings]
{Finitely presented and coherent\\
ordered modules and rings}

\author{F. Wehrung}
\address{C.N.R.S.\\
D\'epartement de Math\'ematiques\\
Universit\'e de Caen\\
14032 Caen Cedex\\
France}
\email{wehrung@math.unicaen.fr}
\urladdr{http://www.math.unicaen.fr/\~{}wehrung}


\begin{abstract}
We extend the usual definition of \emph{coherence}, for modules over
rings, to partially ordered right modules over a large class of
partially ordered rings, called \por s. In this situation,
coherence is equivalent to saying that solution sets of finite systems
of inequalities are finitely generated semimodules.
Coherence for ordered rings and modules, which we call
\po-coherence, has the following features:
\begin{enumerate}
\item Every subring of $\QQ$, and every totally ordered division
ring, is \po-coherent.

\item For a partially ordered right module $\bA$ over a \po-coherent
\por\ $\bR$, $\bA$ is \po-coherent if and only if $\bA$ is a
finitely presented $\bR$-module and $\bA^+$ is a finitely
generated $\bR^+$-semimodule.

\item Every finitely \po-presented partially ordered right module over a
right \po-coherent \por\ is \po-coherent.

\item Every finitely presented \alo\ group is \po-coherent.
\end{enumerate}
\end{abstract}

\maketitle

\section*{Introduction}

Let $\bA$ be a partially ordered right module over a partially
ordered ring $\bR$. We say that $\bA$ is
\emph{\po-coherent}, if the solution sets of finite systems of
inequalities over $\bA$ with unknowns in $\bR$ are finitely generated
\emph{semimodules} over $\bR^+$, that is, monoids closed under
multiplication by scalars in $\bR^+$. This extends the classical
definition of coherence for modules.

In order to avoid pathologies and unwieldy statements, we restrict our
attention to the quite general class of partially ordered right
modules over \emph{\por s}, see Definition~\ref{D:NorOrdRing}. In
particular, every totally ordered ring is a \por. It turns
out that for those modules, the classical notion of \emph{finitely
presented} structure (see Definition~\ref{D:FinPres}), that we shall denote
here by \emph{finitely \po-presented} to emphasize the presence of the
partial ordering, can
conveniently be replaced by the slightly stronger notion of
\emph{finitely related} structure (see Definition~\ref{D:FRPCoh}). Then
we say that a partially ordered right module $\bA$ over a \por\ $\bR$ is
\emph{\po-coherent} (see Definition~\ref{D:FRPCoh}), if every finitely
generated submodule (endowed with the induced ordering) is finitely
related. Equivalently, solution sets of finite homogeneous systems of
inequalities are finitely generated semimodules
(see Theorem~\ref{T:EquivCoh}).

In the case where $\bR$ is, itself, right \po-coherent (that is,
\po-coherent as a partially ordered right module over itself), the
characterization of finitely \po-presented structures takes a much more
wieldy form:

\begin{tha}
$\bA$ is finitely \po-presented if and only if $\bA$ is finitely related,
if and only if $\bA$ is finitely presented as a $\bR$-module and
$\bA^+$ is finitely generated as a $\bR^+$-semimodule.
\end{tha}

Another particular feature of partially ordered modules over a
\po-coherent \por\ is the following:

\begin{thb}
Let $\bR$ be a \por, let $\bA$ be a \porm{\bR}. If $\bR$
is right \po-coherent and $\bA$ is finitely \po-presented, then $\bA$
is \po-coherent.
\end{thb}

The supply of right \po-coherent \por s is quite large:

\begin{thc}
Every subring of $\QQ$ is \po-coherent.
\end{thc}

The main difficulty encountered in the proof of this
result consists of providing a proof for the ring of integers. However,
this is, basically, well-known, and it follows for example from much
more general results due to Grillet (see \cite{Gril76}) and Effros,
Handelman, and Shen (see \cite{EHS80}) about certain partially ordered
abelian groups called \emph{dimension groups}. A similar result can be
proved for ``dimension vector spaces" over a totally ordered field (or
division ring). A consequence of the proof of this result is that
\emph{every totally ordered division ring is a \po-coherent \por}
(Corollary~\ref{C:GEHScoh2}).

The supply of \po-coherent modules is also quite large:

\begin{thd}
Every finitely presented \alo\ group, viewed as a \poag, is \po-coherent.
\end{thd}

Our proof of this result uses the classical representation of finitely
presented \alo\ groups as groups of piecewise linear functions on
polyhedral cones, see \cite{Keim} for a survey.
\medskip

The proofs of the results of this paper are, basically, easy,
especially when they are expressed with the matrix formalism
introduced in Section~\ref{S:Matrices}. However, the results that
they prove are, apparently, nontrivial, and they provide,
moreover, a convenient common platform for handling a number of
different, but related, results that would, otherwise, require separate
proofs. Examples of such statements can be found in \cite{CaWe}.
Furthermore, as mentioned above, the classical notion of a finitely
presented structure is not the most convenient in our context of
ordered structures, which seems to require this switch to finitely
related structures. By Theorem~\ref{T:FinPresEquiv}, the difference
between finitely \po-presented and finitely related
vanishes anyway for ordered modules over \po-coherent rings, which are the
cases that matter for us.

\section{Basic concepts}\label{S:BasConc}

The objects that we shall consider in this paper are, basically,
rings (always associative and unital), and (always unital) right modules
over those rings. Rings and modules will be denoted by capitalized
boldface roman characters, such as $\bR$, $\bS$, $\bA$, $\bB$, while
matrices will be denoted by capitalized lightface roman characters, such
as $M$, $N$, $P$, $Q$, $X$. Elements of a given module will usually be
denoted by lightface lower case roman characters, such as $a$, $b$, $x$,
while scalars in a given ring will be denoted by lightface lower case
greek letters, such as $\alpha$, $\beta$, $\xi$. Since we
shall consider \emph{right} modules, matrices and module homomorphisms
will be written on the \emph{left}.

We shall put $\NN=\ZZ^+\setminus\{0\}$.
If $\bA$ is a right module over a ring $\bR$ and if $\bS$ is a subset of
$\bA$, then we shall denote by $\MM_{m,n}(\bS)$
the set of all matrices with entries in $\bS$, with $m$ rows and $n$
columns, for all $m$, $n\in\ZZ^+$. Furthermore, we shall put
$\MM_n(\bS)=\MM_{n,n}(\bS)$. Thus, for all $p$, $q$, $r\in\ZZ^+$, the
product operation of matrices takes any pair $\vv<X,Y>$, where
$X\in\MM_{p,q}(\bA)$ and
$Y\in\MM_{q,r}(\bR)$, to a matrix $XY\in\MM_{p,r}(\bA)$.
If $a$ is an element of $\bA$, then we shall often identify the
$1\times 1$ matrix $(a)$ with the element $a$ itself.
\smallskip

Our rings and modules will be \emph{partially ordered}. The
topic of the paper is, in fact, designed to cover the context of
\emph{totally ordered} rings and \emph{partially ordered} right
modules over these rings. All the partially ordered rings we shall work
with will be thus defined as follows:

\begin{definition}\label{D:NorOrdRing}
A \emph{\por} is a ring $\vv<\bR,+,\cdot,0,1>$, endowed with a partial
ordering $\leq$ such that the following conditions hold:
\begin{enumerate}
\item $\vv<\bR,+,\cdot,0,1,\leq>$ is a partially ordered ring, that is,
$x\leq y$ implies that $x+z\leq y+z$
and, if $z\geq0$, $xz\leq yz$ and $zx\leq zy$, for all $x$, $y$,
$z\in{\bR}$.

\item $\bR$ is \emph{directed}, that is, $\bR={\bR}^++(-{\bR}^+)$.

\item $\bR$ satisfies $0\leq 1$.
\end{enumerate}
\end{definition}

Note, in particular, that every \emph{totally ordered} ring is
a \por.

We shall consider \emph{partially ordered right modules} over \por s.
For example, every partially ordered ring $\bR$ is also a partially
ordered right module over itself (and a partially ordered left module
as well).

One can define in a classical fashion \emph{quotients} of partially
ordered right modules by convex submodules, see, for example,
\cite[page 122]{Ribe}.

If $\bB$ is a submodule of $\bA$, then we shall put
 \[
 \Conv\bB=
 \{z\in\bA\mid\exists x,y\in{\bB},\ x\leq z\leq y\}.
 \]
It is easy to see that $\Conv\bB$ is the convex submodule of $\bA$
generated by $\bB$. A \emph{finitely generated convex} submodule of $\bA$
is a submodule of the form $\Conv\bB$, where $\bB$ is a finitely
generated submodule of $\bA$.

Our next section, Section~\ref{S:Matrices}, will be devoted to set a
convenient computational framework, based on matrices, for handling
systems of equations and inequalities over partially ordered right
modules. This is done, mainly, in order to avoid carrying everywhere
large amounts of indices which would make the notations considerably
heavier.

\section{Matricial representations of systems of equations
and inequalities}\label{S:Matrices}

We shall fix in this section a \por\ $\bR$ and a
\porm{\bR} $\bA$. In the whole paper, we shall basically consider
\emph{homogeneous} systems of
\emph{equations} and \emph{inequalities}. The typical form of a system of
inequalities is the following:
 \begin{equation}\label{Eq:TypSys}
 \begin{cases}
 a_{11}\xi_1+\cdots+a_{1n}\xi_n&\geq0\\
 a_{21}\xi_1+\cdots+a_{2n}\xi_n&\geq0\\
 \qquad\qquad\vdots\qquad\qquad\vdots&\vdots\\
 a_{m1}\xi_1+\cdots+a_{mn}\xi_n&\geq0,
 \end{cases}
 \end{equation}
where the elements $a_{ij}$ belong to $\bA$ and the elements $\xi_j$
belong to $\bR$. Note that the element
$a_{i1}\xi_1+\cdots+a_{in}\xi_n$ belongs to $\bA$, for all $i\in\fs(m)$.

Now we define matrices $M\in\MM_{m,n}(\bA)$ and $X\in\MM_{n,1}(\bR)$ as
follows:
 \[
 M=\begin{pmatrix}
 a_{11}&a_{12}&\ldots&a_{1n}\\
 a_{21}&a_{22}&\ldots&a_{2n}\\
 \vdots&\vdots&\ddots&\vdots\\
 a_{m1}&a_{m2}&\ldots&a_{mn}
   \end{pmatrix},\qquad\text{and}\qquad
 X=\begin{pmatrix}
 \xi_1\\
 \xi_2\\
 \vdots\\
 \xi_n
   \end{pmatrix}.
 \]
Then the system \eqref{Eq:TypSys} can be written under the following
simple form:
 \begin{equation}\label{Eq:MatSys}
 MX\geq0.
 \end{equation}

If we think of $M$ as \emph{given}, and $X$ as representing the matrix
of \emph{unknowns} of the system \eqref{Eq:TypSys}, then it is natural
that we wish to define the notion of a \emph{general solution} of
\eqref{Eq:TypSys}, and thus of \eqref{Eq:MatSys}. Denote by $\SS$ the
\emph{solution set} of \eqref{Eq:MatSys}, that is,
 \[
 \SS=\{X\in\MM_{n,1}(\bR)\mid MX\geq0\}.
 \]
Then $\SS$ is a $\bR^+$-subsemimodule of $\MM_{n,1}(\bR)$, that is, an
additive submonoid closed under scalar multiplication by elements of
$\bR^+$. Similar considerations hold for systems of \emph{equations}. The
relevant information may be recorded as follows:

\begin{lemma}\label{L:FinGen}
Let $n\in\NN$, let $\SS$ be a subset of $\MM_{n,1}(\bR)$.
Then $\SS$ is a finitely generated right $\bR$-submodule (resp., right
${\bR}^+$-subsemimodule) of
$\MM_{n,1}(\bR)$ if and only if there exist $k\in\NN$ and a matrix
$S\in\MM_{n,k}(\bR)$ such that
 \begin{align*}
 \SS&=\{SY\mid Y\in\MM_{k,1}(\bR)\},\\
 (\text{resp., }\SS&=\{SY\mid Y\in\MM_{k,1}({\bR}^+)\}).\tag*{\qed}
 \end{align*}
\end{lemma}

We shall also consider in this paper systems similar to
 \[
 \begin{cases}
 MX&\geq0\\
 \hfill X&\geq0,
 \end{cases}
 \]
which may be viewed as analogues of \eqref{Eq:MatSys} in the case where
$M$ has entries in $\bR\cup\bA$. We shall call such systems
\emph{mixed systems}.

\section{Finitely related and \po-coherent partially ordered modules}
\label{S:CohMod}

If $\bR$ is a \por\ and if $\bA$ is a \porm{\bR}, we say that a matrix
$M$ with entries in $\bA$ is \emph{spanning}, if the entries of $M$
generate $\bA$ as a $\bR$-submodule. Of particular importance in this
work will be the class of \emph{spanning row matrices}, where a row
matrix is, by definition, a matrix with only one row.

Let us first define, in our context, \emph{finitely \po-presented}
structures:

\begin{definition}\label{D:FinPres}
Let $\bR$ be a \por, let $\bA$ be a \porm{\bR}.
\begin{enumerate}
\item Let $n\in\NN$, let $U\in\MM_{1,n}(\bA)$.
We say that $\bA$ is \emph{finitely \po-presented at} $U$, if the solution
set in $\MM_{n,1}(\bR)$ of the system $UX\geq0$ is a finitely generated
right ${\bR}^+$-subsemimodule of $\MM_{n,1}(\bR)$.

\item $\bA$ is \emph{finitely \po-presented}, if $\bA$ is finitely
\po-presented at \emph{some} spanning row matrix $U$ of elements
of $\bA$ (in particular, $\bA$ is a finitely generated right
$\bR$-module).
\end{enumerate}
\end{definition}

In our context of \por s, this definition agrees with the classical
definition of finitely presented structures in terms of generators and
relations; see, below, the remark (i) following
Definition~\ref{D:FRPCoh}.

More relevant to us will be the following stronger definition of a
\emph{finitely related} structure:

\begin{definition}\label{D:FRPCoh}
Let $\bR$ be a \por, let $\bA$ be a \porm{\bR}.
\begin{enumerate}
\item Let $n\in\NN$, let $U\in\MM_{1,n}(\bA)$.
We say that $\bA$ is \emph{finitely related at} $U$, if the solution set
in $\MM_{n,1}(\bR)$ of the mixed system
 \[
 \begin{cases}
 UX&\geq0\\
 \hfill X&\geq0
 \end{cases}
 \]
is a finitely generated right ${\bR}^+$-subsemimodule of $\MM_{n,1}(\bR)$.

\item $\bA$ is \emph{finitely related}, if $\bA$ is a finitely generated
right $\bR$-module, and $\bA$ is finitely related at
\emph{every} spanning row matrix $U$ of elements of $\bA$.

\item $\bA$ is \emph{\po-coherent}, if every finitely generated
submodule of $\bA$ is finitely
related; that is, if $\bA$ is finitely related at every row matrix of
elements of
$\bA$.
\end{enumerate}

Say that $\bR$ is a \emph{right \po-coherent} \por\ if $\bR$, viewed as a
\porm{\bR}, is \po-coherent.
\end{definition}

Of course, one can define similarly \po-coherent \emph{left} modules over
a \por, and \emph{left \po-coherent} \por s. A \por\ is \emph{\po-coherent},
if it is both right and left \po-coherent.

Several remarks are in order.

\begin{enumerate}
\item Definition~\ref{D:FinPres}
agrees with the classical definition given in universal algebra, stated,
for example, in \cite{Malc} (page 32 for the general definition of an
algebraic system, and page 223 for the definition of a finitely presented
algebraic system). The verification of this easy exercise
strongly uses the fact that $\bR$ is \emph{directed}. For example, this is
already used at the simplest level, to verify that the solution set of the
inequality $(0)\cdot(x)\geq0$ should be a finitely generated right
\emph{subsemimodule} (and not just submodule) of $\bR$. This holds, in
particular, if $\bR$ is directed.

\item The terminology \emph{coherent} is also inspired from classical
module theory:

\begin{definition}\label{D:CohMod}
Let $\bR$ be a ring.
\begin{enumerate}
\item Let $\bA$ be a right $\bR$-module. Say that $\bA$ is \emph{coherent},
if every finitely generated submodule of $\bA$ is finitely presented.

\item Say that $\bR$ is a \emph{right coherent} ring, if $\bR$, viewed as
a right $\bR$-module, is coherent.
\end{enumerate}
\end{definition}

This terminology differs slightly from the one used in
\cite{Sten}, where it is supposed, in addition, that $\bA$ is
itself finitely presented. However, for rings, both terminologies are
equivalent.

\item It may seem slightly awkward to introduce the different but very
similar definitions of finitely \po-presented and finitely related
\porm{\bR}. However, we shall see that in all the cases that will
matter to us, these definitions are equivalent (see
Theorem~\ref{T:FinPresEquiv}).
\end{enumerate}

The general connection between finitely \po-presented and finitely related
\porm{\bR}s is the following:

\begin{proposition}\label{P:FRimplFP}
Let $\bR$ be a \por. Then every finitely related \porm{\bR} is finitely
\po-presented.
\end{proposition}

\begin{proof}
Let $\bA$ be a finitely related \porm{\bR}. Let $n\in\NN$, let
$U\in\MM_{1,n}(\bA)$ be a spanning row matrix of $\bA$. Since $\bA$ is
finitely related and the matrix
$\begin{pmatrix}U& -U\end{pmatrix}$ is a spanning row matrix of $\bA$,
there exist $k\in\NN$ and $P$,
$Q\in\MM_{n,k}(\bR^+)$ such that the matrices $Y$,
$Z\in\MM_{n,1}(\bR^+)$ which satisfy $UY-UZ\geq0$ are exactly the matrices
of the form $Y=PT$, $Z=QT$ for some $T\in\MM_{k,1}(\bR^+)$. Now,
consider an element $X\in\MM_{n,1}(\bR)$. Since $\bR$ is directed, there
are $Y$, $Z\in\MM_{n,1}(\bR^+)$ such that $X=Y-Z$; it follows
that $UX\geq0$ if and only if $UY-UZ\geq0$. Thus, the matrices
$X\in\MM_{n,1}(\bR)$ such that $UX\geq0$ are exactly the matrices of
the form $(P-Q)T$, where $T\in\MM_{k,1}(\bR^+)$.
\end{proof}

\begin{remark}
For all the cases that will matter to us, the converse of
Proposition~\ref{P:FRimplFP} will hold; see
Theorem~\ref{T:FinPresEquiv}.
\end{remark}

By using the formalism of Section~\ref{S:Matrices}, we can easily
recover, in a self-contained fashion, the well-known result that finite
presentability as defined above does not depend on the generating
subset, see V.11, Corollary~7, p. 223 in \cite{Malc}:

\begin{proposition}\label{P:FinPres}
Let $\bA$ be a finitely generated \porm{\bR}. Then $\bA$ is finitely
\po-presented if and only if $\bA$ is finitely \po-presented at
\emph{every} spanning row matrix of $\bA$.
\end{proposition}

\begin{proof}
We prove the nontrivial direction.
So let $U$ and $V$ be spanning row matrices of $\bA$, with, say,
$U\in\MM_{1,m}(\bA)$ and $V\in\MM_{1,n}(\bA)$. Since each of the matrices
$U$ and $V$ has its entries lying in the $\bR$-submodule generated by the
entries of the other, there are matrices $M\in\MM_{m,n}(\bR)$ and
$N\in\MM_{n,m}(\bR)$ such that $V=UM$ and $U=VN$. Suppose that $\bA$ is
finitely \po-presented at $U$, that is, there are
$k\in\NN$ and $S\in\MM_{m,k}(\bR)$ such that
 \begin{equation}\label{Eq:Sol1}
 \{X\in\MM_{m,1}(\bR)\mid UX\geq0\}=
 \{SZ\mid Z\in\MM_{k,1}({\bR}^+)\}.
 \end{equation}
Now let $Y\in\MM_{n,1}(\bR)$ such that $VY\geq0$. This can be written
$UMY\geq0$, thus, by \eqref{Eq:Sol1}, there exists
$Z\in\MM_{k,1}({\bR}^+)$ such that $MY=SZ$. Thus $NMY=NSZ$. Furthermore,
since $\bR$ is directed, there are $Y_0$, $Y_1\in\MM_{n,1}({\bR}^+)$ such
that $Y=Y_0-Y_1$. Therefore, the equality
 \begin{align}
 Y&=(Y-NMY)+NMY\notag\\
  &=(I_n-NM)Y_0+(NM-I_n)Y_1+NSZ\label{Eq:mixed}
 \end{align}
holds, where $I_n$ denotes the identity matrix of order $n$.
Conversely, the equalities $V(I_n-NM)=V-VNM=0$ and $VNS=US\geq0$ hold,
thus every matrix $Y$ of the form \eqref{Eq:mixed} satisfies $VY\geq0$. It
follows that
 \begin{multline*}
 \{Y\in\MM_{n,1}(\bR)\mid VY\geq0\}=\\
 \{(I_n-NM)Y_0+(NM-I_n)Y_1+NSZ\mid\\
 Y_0,\,Y_1\in\MM_{n,1}({\bR}^+)\text{ and }Z\in\MM_{k,1}({\bR}^+)\},
 \end{multline*}
so that $\bA$ is also finitely \po-presented at $V$.
\end{proof}

The following result gives a characterization of coherence for
partially ordered modules:

\begin{theorem}\label{T:EquivCoh}
Let $\bR$ be a \por, let $\bA$ be a \porm{\bR}. Then $\bA$ is
\po-coherent if and only if for all $k$, $m$, $n\in\ZZ^+$ and all
$M\in\MM_{k,m}(\bA)$, $N\in\MM_{k,n}(\bA)$, the solution set of the
following mixed system
\begin{equation}\label{Eq:halfsys}
 \begin{cases}
 MX+NY&\geq 0\\
 \hfill X&\geq 0,
 \end{cases}
\end{equation}
with unknowns $X\in\MM_{m,1}(\bR)$ and $Y\in\MM_{n,1}(\bR)$, is a finitely
generated right $\bR^+$-semimodule.
\end{theorem}

\begin{proof}
We prove the nontrivial direction. So, suppose that $\bA$ is \po-coherent.
We shall first prove a claim:

\begin{sclaim}
For all $m$, $n\in\NN$ and all
$M\in\MM_{m,n}(\bA)$, the solution set
in $\MM_{n,1}(\bR)$ of the mixed system
$\begin{cases}MX\geq0;\quad X\geq 0\end{cases}$ is a finitely generated
right ${\bR}^+$-subsemimodule of $\MM_{n,1}(\bR)$.
\end{sclaim}

\begin{scproof}
We argue by induction on $m$. For $m=1$, the conclusion is just the
assumption on $\bA$. Suppose the claim is true for $m$. Let
$M\in\MM_{m+1,n}(\bA)$, and consider the mixed system
\begin{equation}\label{Eq:samesys}
 \begin{cases}
 MX&\geq 0\\
 \hfill X&\geq 0.
 \end{cases}
\end{equation}
Decompose $M$ into blocks, as
$M=\begin{pmatrix}N\\ U\end{pmatrix}$, where
$N\in\MM_{m,n}(\bA)$ and $U\in\MM_{1,n}(\bA)$. In particular, the mixed
system \eqref{Eq:samesys} is equivalent to the following mixed system:
 \begin{equation}\label{Eq:sys2}
 \begin{cases}
 NX&\geq 0\\
 \hfill UX&\geq0\\
 \hfill X&\geq 0.
 \end{cases}
 \end{equation}
By assumption, the solution set of the mixed system consisting of the last
two subsystems of \eqref{Eq:sys2} above is a finitely generated right
${\bR}^+$-subsemimodule of $\MM_{n,1}(\bR)$, that is, there exist $p\in\NN$
and $S\in\MM_{n,p}(\bR)$ such that
 \begin{equation}\label{Eq:eq1}
 \{X\in\MM_{n,1}({\bR}^+)\mid UX\geq0\}=
 \{SY\mid Y\in\MM_{p,1}({\bR}^+)\}.
 \end{equation}
Note, in particular, that $S\geq0$, and that $NS$ belongs to
$\MM_{m,p}(\bA)$. By induction hypothesis, there are $q\in\NN$ and
$P\in\MM_{p,q}(\bR)$ such that
 \begin{equation}\label{Eq:eq2}
 \{Y\in\MM_{p,1}({\bR}^+)\mid NSY\geq0\}=
 \{PZ\mid Z\in\MM_{q,1}({\bR}^+)\}.
 \end{equation}
Note, in particular, that $P\geq0$.
It follows then easily from \eqref{Eq:eq1} and \eqref{Eq:eq2} that the
solution set of the mixed system \eqref{Eq:sys2} in $\MM_{n,1}(\bR)$ is
equal to
 \[
 \{SPZ\mid Z\in\MM_{q,1}({\bR}^+)\},
 \]
which establishes the induction step.
\end{scproof}

Now we can conclude the proof of Theorem~\ref{T:EquivCoh}.
We consider the mixed system
\begin{equation}\label{Eq:thirdsys}
 \begin{cases}
 MX+NY_0-NY_1&\geq 0\\
 \hfill
  \begin{pmatrix}
  X\\ Y_0\\ Y_1
  \end{pmatrix}
 &\geq 0,
 \end{cases}
\end{equation}
with unknowns $X\in\MM_{m,1}(\bR)$ and
$Y_0$, $Y_1\in\MM_{n,1}(\bR)$. By the Claim, there are $l\in\NN$ and
matrices $P\in\MM_{m,l}(\bR^+)$, $Q_0$, $Q_1\in\MM_{n,l}(\bR^+)$ such that
the solution set of \eqref{Eq:thirdsys} consists of all matrices of the
form
 \[
 \begin{pmatrix}
 PT\\ Q_0T\\ Q_1T
 \end{pmatrix},
 \]
for $T\in\MM_{l,1}(\bR^+)$. But every $Y\in\MM_{n,1}(\bR)$ can be written
as a difference $Y=Y_0-Y_1$ where $Y_0$, $Y_1\in\MM_{n,1}(\bR^+)$. Hence
the solution set of \eqref{Eq:halfsys} consists of all matrices of the
form
 \[
 \begin{pmatrix}
 PT\\ (Q_0-Q_1)T
 \end{pmatrix},
 \]
for $T\in\MM_{l,1}(\bR^+)$.
\end{proof}

\begin{corollary}\label{C:poCohCoh}
Let $\bR$ be a \por, let $\bA$ be a \porm{\bR}. If $\bA$ is \po-coherent,
then $\bA$ is coherent.
\end{corollary}

\begin{proof}
Let $n\in\NN$ , let $U\in\MM_{1,n}(\bA)$. We must prove that the set
 \[
 \SS=\{X\in\MM_{n,1}(\bR)\mid UX=0\}
 \]
is a finitely generated right $\bR$-module. Since the system $UX=0$ is
equivalent to the system
 \[
 \begin{pmatrix}
 U\\ -U
 \end{pmatrix}
 X\geq0,
 \]
this follows immediately from Theorem~\ref{T:EquivCoh}.
\end{proof}

\section{Right \po-coherent \por s}

This section is devoted to the study of a few characterizations of
\po-coherent partially ordered right modules. We start with the following
result.

\begin{theorem}\label{T:SubmFP}
Let $\bR$ be a \por, let $\bA$ be a \porm{\bR}. If $\bR$
is right \po-coherent and $\bA$ is finitely \po-presented, then $\bA$
is \po-coherent.
\end{theorem}

\begin{proof}
We first establish two claims:

\setcounter{claim}{0}
\begin{claim}
$\bA$ is finitely related.
\end{claim}

\begin{cproof}
By assumption, $\bA$ is finitely generated. Let $U$ be a spanning row matrix of $\bA$, say, $U\in\MM_{1,n}(\bA)$, where $n\in\NN$. By Proposition~\ref{P:FinPres}, $\bA$ is finitely \po-presented at $U$, that is, there are $p\in\NN$ and $S\in\MM_{n,p}(\bR)$ such that
 \begin{equation}\label{Eq:SolPR}
 \{X\in\MM_{n,1}(\bR)\mid UX\geq0\}=\{SY\mid Y\in\MM_{p,1}({\bR}^+)\}.
 \end{equation}
Note, in particular, that $S$ is a matrix with entries in $\bR$. By
Theorem~\ref{T:EquivCoh}, since $\bR$ is right \po-coherent, there exist
$q\in\NN$ and $T\in\MM_{p,q}({\bR}^+)$ such that
 \begin{equation}\label{Eq:inporing}
 \{Y\in\MM_{p,1}({\bR}^+)\mid SY\geq0\}
 =\{TZ\mid Z\in\MM_{q,1}({\bR}^+)\}.
 \end{equation}
Note, in particular, that $T\geq0$ and $ST\geq0$.
It follows easily from \eqref{Eq:SolPR} and \eqref{Eq:inporing} that
the equality
 \[
 \{X\in\MM_{n,1}({\bR}^+)\mid UX\geq0\}=
 \{STZ\mid Z\in\MM_{q,1}({\bR}^+)\}
 \]
holds. Therefore, $\bA$ is finitely related at $U$.
\end{cproof}

\begin{claim}
Every finitely generated submodule of $\bA$ is finitely \po-pre\-sent\-ed.
\end{claim}

\begin{cproof}
Let $\bB$ be a finitely generated submodule of $\bA$. Let $U$ (resp., $V$)
be a spanning row matrix of $\bA$ (resp., $\bB$), say, $U\in\MM_{1,m}(\bA)$
and $V\in\MM_{1,n}(\bB)$. Since $\bB$ is contained in $\bA$, there exists
a matrix $M\in\MM_{m,n}(\bR)$ such that $V=UM$. We define sets $\SS$ and
$\TT$ of matrices as follows:
 \[
 \SS=\{X\in\MM_{m,1}(\bR)\mid UX\geq0\},\quad
 \TT=\{Y\in\MM_{n,1}(\bR)\mid VY\geq0\}.
 \]
By assumption on $\bA$, $\SS$ is a finitely generated
${\bR}^+$-subsemimodule of $\MM_{m,1}(\bR)$, that is, there exist $k\in\NN$
and $S\in\MM_{m,k}(\bR)$ such that
$\SS=\{SZ\mid Z\in\MM_{k,1}({\bR}^+)\}$.
Now consider the following mixed system:
 \begin{equation}\label{Eq:MYSZ}
 \begin{cases}
 MY-SZ&=0\\
 \hfill Z&\geq0
 \end{cases}
 \end{equation}
with unknowns $Y$ in $\MM_{n,1}(\bR)$ and $Z$ in $\MM_{k,1}(\bR)$.
Since $\bR$ is right \po-coherent, there are, by Theorem~\ref{T:EquivCoh},
$l\in\NN$ and matrices $P\in\MM_{n,l}(\bR)$, $Q\in\MM_{k,l}(\bR)$ such that
the solution set of \eqref{Eq:MYSZ} equals the set
 \[
 \left\{
  \begin{pmatrix}
  P\\
  Q
  \end{pmatrix}
  T\mid T\in\MM_{l,1}({\bR}^+)
 \right\}.
 \]
But a matrix $Y\in\MM_{n,1}(\bR)$ belongs to $\TT$ if and only if
$VY\geq0$, that is, $UMY\geq0$, which can be written $MY\in\SS$, or:
there is a matrix $Z$ such that $Y$ and $Z$ satisfy \eqref{Eq:MYSZ}.
It follows that $\TT=\{PT\mid T\in\MM_{l,1}({\bR}^+)\}$. This proves that
$\TT$ is a finitely generated ${\bR}^+$-subsemimodule of $\MM_{n,1}(\bR)$.
\end{cproof}

We can now conclude the proof of Theorem~\ref{T:SubmFP}. Let $\bB$ be a
finitely generated submodule of $\bA$. By Claim~2, $\bB$ is finitely
\po-presented. By Claim~1, applied to $\bB$, $\bB$ is finitely related.
\end{proof}

We can now give a characterization of finite \po-presentability that
separates the roles of the algebraic structure and of the ordering:

\begin{theorem}\label{T:FinPresEquiv}
Let $\bR$ be a right \po-coherent \por. For every
\porm{\bR} $\bA$, the following are equivalent:
\begin{enumerate}
\item $\bA$ is finitely \po-presented.

\item $\bA$ is finitely related.

\item $\bA$ is a finitely presented right $\bR$-module, and
${\bA}^+$ is a finitely generated right ${\bR}^+$-semimodule.
\end{enumerate}
\end{theorem}

\begin{proof}
(i)$\Rightarrow$(ii) follows immediately from Theorem~\ref{T:SubmFP}.

(ii)$\Rightarrow$(i) has been proved in Proposition~\ref{P:FRimplFP}.

(i)$\Rightarrow$(iii) Suppose that (i) holds. By Theorem~\ref{T:SubmFP},
$\bA$ is \po-coherent. By Corollary~\ref{C:poCohCoh}, $\bA$ is coherent.
Since $\bA$ is finitely \po-presented, it is finitely generated. It
follows that $\bA$ is finitely presented.

It remains to prove that $\bA^+$ is finitely generated over $\bR^+$. By
assumption, there are $n\in\NN$ and a spanning row
matrix $U\in\MM_{1,n}(\bA)$ at which $\bA$ is finitely \po-presented.
Hence, there exist $k\in\NN$ and $S\in\MM_{n,k}(\bR)$ such that
 \begin{equation}\label{Eq:PosCone}
 \{X\in\MM_{n,1}(\bR)\mid UX\geq0\}=\{SY\mid Y\in\MM_{k,1}({\bR}^+)\}.
 \end{equation}
Since $U$ is a spanning row matrix of $\bA$, every element of $\bA$ has
the form $UX$ where $X\in\MM_{n,1}(\bR)$. Hence, by \eqref{Eq:PosCone},
every element of ${\bA}^+$ has the form $USY$, where
$Y\in\MM_{k,1}(\bR^+)$; thus it is a positive linear combination of the
finite set $T=\{USE_l\mid 1\leq l\leq k\}$, where
$\vv<E_l\mid 1\leq l\leq k>$ denotes the canonical basis of
$\MM_{k,1}(\bR)$. Since $0\leq1$ in $\bR$, the elements $USE_l$, for
$1\leq l\leq k$, all belong to ${\bA}^+$; thus $T$ is a finite generating
subset of ${\bA}^+$.

(iii)$\Rightarrow$(i) Suppose that (iii) holds. Let $m$, $n\in\NN$ such
that there exists a finite generating subset $\{a_1,\ldots,a_m\}$ of
$\bA$ as a right $\bR$-module, and a finite generating subset
$\{b_1,\ldots,b_n\}$ of ${\bA}^+$ as a right ${\bR}^+$-semimodule.
We define row matrices $U$ and $V$ by
 \[
 U=\begin{pmatrix}
   a_1&\cdots&a_m
   \end{pmatrix}
 \qquad\text{and}\quad
 V=\begin{pmatrix}
   b_1&\cdots&b_n
   \end{pmatrix}.
 \]
We can assume that $\bA$ is finitely presented at $U$.
Since the entries of $U$ generate $\bA$ as a right $\bR$-module, there
exists $M\in\MM_{m,n}(\bR)$ such that $V=UM$. Since $\bA$ is finitely
presented at $U$, the solution set $\SS$ in
$\MM_{m,1}(\bR)$ of the matrix equation $UX=0$ is a finitely generated
right $\bR$-submodule of $\MM_{m,1}(\bR)$, thus, since $\bR$ is directed,
it is also a finitely generated right $\bR^+$-subsemimodule of
$\MM_{m,1}(\bR)$.

Put $\mathcal{P}=\{X\in\MM_{m,1}(\bR)\mid UX\geq0\}$. To conclude the
proof, it suffices to prove that $\mathcal{P}$ is a finitely generated
right ${\bR}^+$-subsemimodule of $\MM_{m,1}(\bR)$. Since $\SS$ is a
finitely generated right ${\bR}^+$-subsemimodule of $\MM_{m,1}(\bR)$, it
suffices to establish the following formula:
 \begin{equation}\label{Eq:PSplus}
 \mathcal{P}=\{MY\mid Y\in\MM_{n,1}({\bR}^+)\}+\SS.
 \end{equation}
The containment from the right into the left follows from the
fact that $UM=V\geq0$. Conversely, let $X\in\MM_{m,1}(\bR)$ such that
$UX\geq0$. Since the entries of $V$ generate ${\bA}^+$ as a right
${\bR}^+$-subsemimodule, there exists $Y\in\MM_{n,1}({\bR}^+)$ such that
$UX=VY$. This can be written $UX=UMY$, so that $X-MY\in\SS$. This
completes the proof of \eqref{Eq:PSplus}.
\end{proof}

\begin{remark}
Unlike (i)$\Rightarrow$(ii), several implications in the proof of
Theorem~\ref{T:FinPresEquiv} do not use the full strength of the
hypothesis that $\bR$ is a right \po-coherent \por. For example,
(ii)$\Rightarrow$(i) and (iii)$\Rightarrow$(i)
hold for arbitrary \por s.
\end{remark}

\section{Coherence preservation properties for ordered modules
and rings}\label{S:Clos}

\begin{lemma}\label{L:FinIntSem}
Let $\bR$ be a \por, let $\bA$ be a \porm{\bR}. If $\bA$ is
\po-coherent, then the set of all finitely generated $\bR$-modules
(resp., ${\bR}^+$-subsemimodules) of $\bA$ is closed under finite
intersection.
\end{lemma}

\begin{remark}
The part about $\bR$-modules is not really a property of partially
ordered modules, and it holds, in fact, in the case where $\bA$ is a
coherent module over a ring $\bR$.
\end{remark}

\begin{proof}
We give a proof for the semimodule part; the proof for modules is
similar. Let $\bB$ and $\bC$ be finitely generated
${\bR}^+$-subsemimodules of $\bA$. Let $U$ (resp., $V$) be a spanning row
matrix of $\bB$ (resp.,
$\bC$), say, $U\in\MM_{1,m}(\bB)$ and $V\in\MM_{1,n}(\bC)$. Since $\bA$
is \po-coherent, the solution set of the following mixed system
 \[
 \begin{cases}
 UX-VY&=0\\
 \hfill X&\geq0\\
 \hfill Y&\geq0,
 \end{cases}
 \]
with unknowns $X\in\MM_{m,1}(\bR)$ and $Y\in\MM_{n,1}(\bR)$,
is, by Theorem~\ref{T:EquivCoh}, a finitely generated
${\bR}^+$-subsemimodule of $\MM_{m,1}(\bR)\times\MM_{n,1}(\bR)$, that is,
there exist $k\in\NN$ and matrices $P\in\MM_{m,k}({\bR}^+)$,
$Q\in\MM_{n,k}({\bR}^+)$ such that
 \[
  \left\{
   \begin{pmatrix}
   X\\
   Y
   \end{pmatrix}
   \in\MM_{m+n,1}({\bR}^+)\mid UX=VY
  \right\}
  =
  \left\{
   \begin{pmatrix}
   PZ\\
   QZ
   \end{pmatrix}
   \mid Z\in\MM_{k,1}({\bR}^+)
  \right\}.
 \]
In particular, $UP=VQ$. If we denote this matrix by $S$, it is then
clear that
 \begin{equation*}
 \bB\cap\bC=\{SZ\mid Z\in\MM_{k,1}({\bR}^+)\}.\tag*{\qed}
 \end{equation*}
\renewcommand{\qed}{}
\end{proof}

\begin{lemma}\label{L:FinDimCoh}
Let $\bR$ be a right \po-coherent \por. Then
$\bR^n$, partially ordered coordinatewise, is a \po-coherent \porm{\bR},
for all $n\in\ZZ^+$.
\end{lemma}

\begin{proof}
This follows immediately from the definition of \po-coherence, and
Theorem~\ref{T:EquivCoh}.
\end{proof}

We deduce the following preservation property:

\begin{proposition}\label{P:nRRcoh}
Let $\bR$ be a right \po-coherent \por. Then every direct sum of
\po-coherent \porm{\bR}s, partially ordered coordinatewise, is
\po-coherent.
\end{proposition}

\begin{proof}
Let $\vv<\bA_i\mid i\in I>$ be a family of \po-coherent \porm{\bR}s; put
$\bA=\bigoplus_{i\in I}\bA_i$. Let $n\in\NN$, let $a_1$, \dots,
$a_n$ be elements of $\bA$. Let $J$ be the union of the supports of all
the $a_k$; so, $J$ is a finite subset of $I$.
We write $a_k=\vv<a_{i,k}\mid i\in I>$, for all $k\in\fs(n)$; in
particular, $a_{i,k}=0$ if $i\in I\setminus J$. It follows that for all
$\xi_1$, \dots, $\xi_n\in{\bR}^+$, $\sum_{k=1}^na_k\xi_k\geq0$ if and only
if $\sum_{k=1}^na_{i,k}\xi_k\geq0$, for all
$i\in J$. Since $\bA_i$ is \po-coherent, the
set of all $\vv<\xi_1,\ldots,\xi_n>\in({\bR}^+)^n$ satisfying the
inequality $\sum_{k=1}^na_{i,k}\xi_k\geq0$ is a finitely generated
${\bR}^+$-subsemimodule of $({\bR}^+)^n$. We conclude the proof by Lemmas
\ref{L:FinIntSem} and \ref{L:FinDimCoh}.
\end{proof}

\begin{corollary}\label{C:nRRcoh}
Let $\bR$ be a right \po-coherent \por. Then every free right $\bR$-module,
endowed with the componentwise ordering, is a \po-coherent \porm{\bR}.\qed
\end{corollary}

It is obvious that every submodule of a \po-coherent partially ordered
module is \po-coherent. The situation is slightly more complicated for 
quotients by convex submodules. In particular, we prove in
Example~\ref{Ex:lgroups} that the quotient of a \po-coherent
\poag\ by a convex, directed subgroup may not be \po-coherent.

We shall now see that this problem
does not arise for quotients by \emph{finitely generated} convex
submodules.

\begin{proposition}\label{P:ConvQuot}
Let $\bR$ be a \por, let $\bA$ be a \porm{\bR}, let
$\bF$ be a finitely generated convex submodule of $\bA$.
If $\bA$ is \po-coherent, then $\bA/\bF$ is \po-coherent.
\end{proposition}

\begin{proof}
Let $\bB$ be a finitely generated submodule of $\bA$ such that
$\bF=\Conv\bB$ (see Section~\ref{S:BasConc}). Furthermore,
denote by $\pi$ the natural projection from $\bA$ onto $\bA/\bF$. Let
$V$ be a spanning row matrix of $\bB$, say, $V\in\MM_{1,n}(\bB)$,
for some $n\in\NN$.

Let $\overline{U}$ be a row matrix in $\bA/\bF$, say,
$\overline{U}\in\MM_{1,m}(\bA/\bF)$, for some $m\in\NN$. There exists
$U\in\MM_{1,m}(\bA)$ such that $\overline{U}=\pi U$.
We must prove that the solution set in $\MM_{m,1}(\bR)$ of the mixed system
$\begin{cases}\overline{U}X\geq0;\quad X\geq0\end{cases}$,
is a finitely generated ${\bR}^+$-subsemimodule of $\MM_{m,1}(\bR)$.

Consider the following mixed system, with unknowns $X\in\MM_{m,1}(\bR)$ and
$Y\in\MM_{n,1}(\bR)$:
 \begin{equation}\label{Eq:UXVY}
  \begin{cases}
  UX+VY&\geq0\\
  \hfill X&\geq0.
  \end{cases} 
 \end{equation}
Since $\bA$ is \po-coherent, there are,
by Theorem~\ref{T:EquivCoh}, matrices $P\in\MM_{m,k}(\bR)$ and
$Q\in\MM_{n,k}(\bR)$ such that the solutions of \eqref{Eq:UXVY} are
exactly the matrices of the form
$\begin{pmatrix}PZ\\ QZ\end{pmatrix}$, where $Z\in\MM_{k,1}({\bR}^+)$.
But for all $X\in\MM_{m,1}({\bR}^+)$, $\overline{U}X\geq0$ if and only if
there exists $b\in\bF$ such that $UX+b\geq0$, that is, there exists $Y$
such that $X$ and $Y$ satisfy \eqref{Eq:UXVY}. Therefore, the elements
$X$ of $\MM_{m,1}({\bR}^+)$ such that $\overline{U}X\geq0$ are exactly the
matrices of the form $PZ$, where $Z\in\MM_{k,1}({\bR}^+)$.
\end{proof}

We now turn to preservation results for rings:

\begin{proposition}\label{P:MatCoh}
Let $\bR$ be a right \po-coherent \por. Then the matrix
ring $\MM_m(\bR)$, ordered componentwise, is right \po-coherent, for all
$m\in\NN$.
\end{proposition}

\begin{proof}
Put $\bS=\MM_m(\bR)$.
Let $n\in\NN$, let $a_1$, \dots, $a_n\in{\bS}$. For all $\xi_1$,
\dots, $\xi_n\in{\bS}^+$, the condition $\sum_{j=1}^na_j\xi_j\geq0$ can
be expressed by a system of inequalities, with coefficients all the
$nm^2$ entries of all the $a_j$, and unknowns all the $nm^2$ entries of
all the $\xi_j$. Since $\bR$ is right \po-coherent, the solution set of
this system is, by Theorem~\ref{T:EquivCoh}, a finitely generated
${\bR}^+$-subsemimodule of $\MM_{n,1}(\bS)$. Since it is also a
$\bS^+$-subsemimodule of $\MM_{n,1}(\bS)$, it is,
\emph{a fortiori}, a finitely generated
${\bS}^+$-subsemimodule of $\MM_{n,1}(\bS)$.
\end{proof}

We shall make use in Section~\ref{S:SubRat} of the following result.

\begin{proposition}\label{P:ExtCoh}
Let $\bS$ be a \por, let $\bR$ be a directed subring of $\bS$.
We define subsets $\bD$ and $\bU$ of $\bR$ as follows:
 \begin{align*}
 \bD&=\{d\in\bR^+\mid\forall x\in\bS,\ dx\geq0\Rightarrow x\geq0\},\\
 \bU&=\{u\in\bR^+\mid\exists v\in\bS^+,\ uv=1\}.
 \end{align*}

We suppose that for all $s\in\bS$, there exist $d\in\bD$ and $u\in\bU$ such
that $ds\in\bR$ and $su\in\bR$.

Then the following holds:
if $\bR$ is right \po-coherent, then $\bS$ is right \po-coherent.
\end{proposition}

\begin{proof}

We start with a claim:

\begin{sclaim}
Let $n\in\NN$, let $X\in\MM_{1,n}(\bS)$. Then there are $d\in\bD$ and
$u\in\bU$ such that both $dX$ and $Xu$ belong to $\MM_{1,n}(\bR)$.
\end{sclaim}

\begin{scproof}
Easy by induction on $n$, by using the assumption on $\bD$ and $\bU$, and
the obvious fact that both $\bD$ and $\bU$ are closed under multiplication.
\end{scproof}

Let $n\in\NN$ and let $U\in\MM_{1,n}(\bS)$ be a row matrix of
$\bS$. We prove that the solution set in $\MM_{n,1}(\bS)$ of the
following mixed system
 \[
 \begin{cases}
 UX&\geq0\\
 \hfill X&\geq0
 \end{cases}
 \]
is a finitely generated ${\bS}^+$-subsemimodule of $\MM_{n,1}(\bS)$.
By the Claim, there exists $d\in\bD$ such that $dU\in\MM_{1,n}(\bR)$.
Note that for all $X\in\MM_{n,1}(\bR^+)$, $UX\geq0$ if and only if
$dUX\geq0$. Hence, we may assume without loss of generality that
$U$ belongs to $\MM_{1,n}(\bR)$.
Since $\bR$ is right \po-coherent, there exist $k\in\NN$ and
$P\in\MM_{n,k}(\bR)$ such that
 \begin{equation}\label{Eq:UXPYR}
 \{X\in\MM_{n,1}({\bR}^+)\mid UX\geq0\}=
 \{PY\mid Y\in\MM_{k,1}({\bR}^+)\}.
 \end{equation}
Thus, in order to conclude the proof, it suffices to prove that the
following holds:
 \begin{equation}\label{Eq:UXPYS}
 \{X\in\MM_{n,1}({\bS}^+)\mid UX\geq0\}=
 \{PY\mid Y\in\MM_{k,1}({\bS}^+)\}.
 \end{equation}
Since $P\geq0$ and $UP\geq0$, the containment from the right into the left
is obvious. Conversely, let $X\in\MM_{n,1}({\bS}^+)$ such that $UX\geq0$.
By the Claim, there exists $u\in\bU$ such that
$Xu\in\MM_{n,1}({\bR}^+)$. Since $u\geq0$, we also have $UXu\geq0$ and
$Xu\geq0$. By definition, there exists $v\in{\bS}^+$ such that
$uv=1$. By \eqref{Eq:UXPYR}, there exists $Y\in\MM_{k,1}({\bR}^+)$ such
that $Xu=PY$. It follows that $X=Xuv=PYv$. Note that $Yv\geq0$. This
concludes the proof of \eqref{Eq:UXPYS}.
\end{proof}

\section{Subrings of the rationals}\label{S:SubRat}

We shall prove in this section (Theorem~\ref{T:CohSubQ}) that every
subring of $\QQ$ is \po-coherent. The hard core of this fact is
the case of the ring $\ZZ$ of integers. The result follows then from a
much more general principle, due, in the context of commutative
semigroups, to Grillet, see \cite{Gril76}, and, in the context of
\poag s, to Effros, Handelman, and Shen, see \cite{EHS80}, and also
\cite{Gpoag}. In order to formulate this result, we need a few
definitions.

\begin{definition}\label{D:DimGroup}
Let $\bG$ be a \poag.
\begin{enumerate}
\item $\bG$ is \emph{directed}, if $\bG=\bG^++(-\bG^+)$.

\item $\bG$ is \emph{unperforated}, if $mx\geq0$ implies that
$x\geq0$, for all $m\in\NN$ and all $x\in\bG$.

\item $\bG$ satisfies the \emph{interpolation property},
or $\bG$ is an \emph{interpolation group}, if for all
$a_0$, $a_1$, $b_0$, $b_1\in\bG$ such that $a_0,a_1\leq b_0,b_1$,
there exists $x\in\bG$ such that $a_0,a_1\leq x\leq b_0,b_1$.

\item $\bG$ is a \emph{dimension group}, if it is a directed, unperforated
interpolation group.
\end{enumerate}
\end{definition}

An important particular class of dimension groups is the class of
\emph{simplicial groups}; a \poag\ is simplicial, if it is isomorphic
to some $\ZZ^n$, with $n\in\ZZ^+$. The Grillet, Effros, Handelman,
and Shen Theorem states that a \poag\ is a dimension group if and only if
it is a direct limit of simplicial groups and positive homomorphisms (that
is, order-preserving group homomorphisms).

\begin{proposition}\label{P:GEHSeqsys}
Let $m$, $n\in\NN$, let $M\in\MM_{m,n}(\ZZ)$. Then there exist
$k\in\NN$ and a matrix $S\in\MM_{k,m}(\ZZ)$ such that for every
unperforated interpolation group $\bG$, the solution set in
$\MM_{1,m}(\bG^+)$ of the system $XM\geq0$ is exactly the set
\begin{equation}
\{YS\mid Y\in\MM_{1,k}(\bG^+)\}.\tag*{\qed}
\end{equation}
\end{proposition}

The case where $m=1$, sufficient to deduce the Grillet, Effros, Handelman,
and Shen Theorem, is done in \cite[Proposition~3.15]{Gpoag}. The general
case is easily deduced by induction, in a similar fashion as in the proof
of Theorem~\ref{T:EquivCoh}.

By applying Proposition~\ref{P:GEHSeqsys} to $\bG=\ZZ$, we obtain the
following result:

\begin{corollary}\label{C:GEHScoh}
The ring $\ZZ$ of integers, endowed with its natural ordering, is
\po-coherent.\qed
\end{corollary}

By using Proposition~\ref{P:ExtCoh}, we thus obtain a large class of
\po-coherent \por s:

\begin{theorem}\label{T:CohSubQ}
Every subring of $\QQ$ is \po-coherent.
\end{theorem}

\begin{proof}
Let $\bR$ be a subring of $\QQ$.
By Corollary~\ref{C:GEHScoh}, the ordered ring $\ZZ$ is \po-coherent. Since
$\ZZ$ is a directed subring of $\bR$, it suffices to prove that the
containment of ordered rings $\ZZ\subseteq\bR$ satisfies the
assumption of Proposition~\ref{P:ExtCoh}. We define subsets
$\bD$ and $\bU$ of $\bR^+$ as in Proposition~\ref{P:ExtCoh}. Then we have
$\bD=\NN$, thus, since $\bR\subseteq\QQ$, (i) is obvious.
Now we prove (ii). So let $x\in\bR$. Write $x=p/q$, where $p\in\ZZ$,
$q\in\NN$ and $p$ and $q$ are coprime. By Bezout's Theorem, there are
integers $u$ and $v$ such that $up+vq=1$. Therefore, the equality
 \[
 \frac{1}{q}=\frac{up+vq}{q}=ux+v\in\bR
 \]
holds, so that $q\in\bU$. Moreover, $qx=p\in\ZZ$, which completes the
verification of (ii).
\end{proof}

\begin{example}\label{E:CohSubR}
We construct a subring of the ring $\RR$ of real numbers that is not
even coherent as a \emph{ring}. Let $u$, $v$ and $x_n$ ($n\in\ZZ^+$)
be algebraically independent elements of $\RR$ over $\QQ$,
and let $\bR$ be the subring of $\RR$ generated by the elements $u$,
$v$, and $ux_n$, $vx_n$ for all $n\in\ZZ^+$. Note, in particular,
that no $x_n$ belongs to $\bR$. Then it is not difficult to verify
that the set $\boldsymbol{I}$ defined by
 \[
 \boldsymbol{I}=\{\vv<x,y>\in\bR\times\bR\mid ux=vy\}
 \]
is the $\bR$-submodule of $\bR\times\bR$ generated by all pairs
$\vv<v,u>$ and $\vv<vx_n,ux_n>$, for $n\in\ZZ^+$, and that it is not
finitely generated. Therefore, $\bR$ is not coherent.
\end{example}

\section{Totally ordered division rings}

In the previous section, we have seen that, in particular, the ordered
field of rationals is a \po-coherent \por.
In this section, we shall see that \emph{every} totally ordered field
is \po-coherent. In fact, the commutativity will not be used in the proof,
so that this result will hold for totally ordered \emph{division
rings}.

We first state the analogue, for totally ordered division rings, of
Definition~\ref{D:DimGroup}:

\begin{definition}\label{D:DimVs}
Let $\bK$ be a totally ordered division ring. A \emph{right \povs}, or,
simply, a \emph{\povs}, over $\bK$ is, by definition, a \porm{\bK}. A
\emph{dimension vector space} over $\bK$ is a directed \povs\ over
$\bK$ satisfying the interpolation property.
\end{definition}

Note that a \povs\ is automatically unperforated, so that every
dimension vector space is also a dimension group. We also extend naturally
to vector spaces the definition, stated in Section~\ref{S:SubRat}, of a
simplicial group: a \emph{simplicial vector space} over $\bK$ is a
\povs\ over $\bK$ which is isomorphic to some $\bK^n$ (endowed with the
natural ordering), for some $n\in\ZZ^+$.

The analogue for \povs s over $\bK$ of the Grillet,
Effros, Handelman, and Shen Theorem is then that for every totally ordered
division ring $\bK$, the dimension vector spaces over $\bK$ are exactly
the direct limits of simplicial vector spaces and positive homomorphisms
of $\bK$-vector spaces.

Similarly, we can formulate for totally ordered division rings the
following analogue of Proposition~\ref{P:GEHSeqsys}.

\begin{proposition}\label{P:GEHSeqsys2}
Let $\bK$ be a totally ordered division ring.
Let $m$, $n\in\NN$, let $M\in\MM_{m,n}(\bK)$. Then there exist
$k\in\NN$ and a matrix $S\in\MM_{k,m}(\bK)$ such that for every
interpolation vector space $\bE$ over $\bK$, the solution set in
$\MM_{1,m}(\bE^+)$ of the system $XM\geq0$ is exactly the set
 \begin{equation}
 \{YS\mid Y\in\MM_{1,k}(\bE^+)\}.\tag*{\qed}
 \end{equation}
\end{proposition}

\begin{proof}
The hard core of the proof consists, similarly as for
Proposition~\ref{P:GEHSeqsys}, to prove that if $\bE$ is an interpolation
vector space over $\bK$, then, for all
$n\in\NN$ and all $p_1$, \dots, $p_n\in\bK$, the solution set in
$(\bE^+)^n$ of the equation
 \begin{equation}\label{Eq:p1nx}
 x_1p_1+\cdots+x_np_n=0
 \end{equation}
is the set of linear combinations with coefficients from $\bE^+$ of
a certain finite subset of $(\bK^+)^n$. Actually, the
proof of this fact is much easier than for groups, because of the possibility of dividing by any nonzero element of $\bK$.

Put $U=\{i\in\fs(n)\mid p_i>0\}$, $V=\{i\in\fs(n)\mid p_i<0\}$,
and $q_i=-p_i$, for all $i\in V$. Note that $U\cap V=\varnothing$.
Then \eqref{Eq:p1nx} can be written as follows:
 \[
 \sum_{i\in U}x_ip_i=\sum_{j\in V}x_jq_j
 \]
Since $\bE^+$ satisfies the refinement property
(see \cite[Proposition~2.2(c)]{Gpoag}), this is equivalent to saying
that there are elements $z_{ij}\in\bE^+$, for $\vv<i,j>\in U\times V$,
such that
 \begin{equation}\label{Eq:piqjz}
 \begin{aligned}
 x_ip_i=\sum_{j\in V}z_{ij},&\qquad\text{for all }i\in U,\\
 x_jq_j=\sum_{i\in U}z_{ij},&\qquad\text{for all }j\in V.
 \end{aligned}
 \end{equation}
Dividing by $p_i$ (resp., $q_j$) yields the following equivalent form
of \eqref{Eq:piqjz}:
 \[
 \begin{aligned}
 x_i=\sum_{j\in V}z_{ij}p_i^{-1},&\qquad\text{for all }i\in U,\\
 x_j=\sum_{i\in U}z_{ij}q_j^{-1},&\qquad\text{for all }j\in V.\\
 \end{aligned}
 \]
Thus the solution set in $(\bE^+)^n$ of \eqref{Eq:p1nx} is the set of
linear combinations with coefficients from $\bE^+$ of a subset of
$(\bK^+)^n$ with $n+|U|\cdot|V|-|U\cup V|$ elements.

For general $m$, the conclusion follows from an easy induction argument, similar to the one used in the proof of Theorem~\ref{T:EquivCoh}.
\end{proof}

\begin{corollary}\label{C:GEHScoh2}
Every totally ordered division ring is \po-coherent.\qed
\end{corollary}

\section{\po-coherence of finitely presented \alo\ groups}

The main goal of this section is to provide a proof for the following
result:

\begin{theorem}\label{T:FinPres}
Every finitely presented \alo\ group, viewed as a \poag, is \po-coherent.
\end{theorem}

We shall need to recall some basic results about \alo\ groups. The
reader may consult the survey article \cite{Keim} for more information.

\subsection{Finitely generated \alo\ groups and polyhedral cones}

We refer to \cite{Bern,Wein1,Wein2} for more information about this
section. Let $n$ be a positive integer. The abelian group $\ZZ^{\ZZ^n}$ of
all maps from $\ZZ^n$ to
$\ZZ$, endowed with its coordinatewise ordering, is an \alo\ group. We
consider the sub-lattice-ordered group $\FA(n)$ of $\ZZ^{\ZZ^n}$ generated
by the $n$ canonical projections $p_i\colon\ZZ^n\to\ZZ$, for
$1\leq i\leq n$. It is
then well-known that $\vv<\FA(n),p_1,\ldots,p_n>$ is the free
\alo\ group on $n$ generators. Every element $f$ of $\FA(n)$ can be
decomposed as
 \begin{equation}\label{Eq:Decomp}
 f=\bigwedge_{1\leq i\leq k}f_i-\bigwedge_{k+1\leq i\leq k+l}f_i,
 \end{equation}
where $k$, $l$ are positive integers and all the $f_i$ are linear
functionals from $\ZZ^n$ to $\ZZ$. In particular, $f$ can be extended to a
unique \emph{positively homogeneous} map from $\QQ^n$ to $\QQ$, defined by
\eqref{Eq:Decomp} where the $f_i$ are replaced by their
unique extensions to linear functionals on $\QQ^n$. Throughout
this section, we shall identify whenever needed an element of $\FA(n)$ with its positively homogeneous extension to $\QQ^n$.

As in \cite{Bake,Keim}, we define a \emph{\cpc} of $\QQ^n$ to be a finite
intersection of closed half-spaces of $\QQ^n$, where, as usual, a closed
half-space is a subset of the form $\{x\mid p(x)\geq0\}$, where $p$ is a
nonzero linear functional on $\QQ^n$. Further, we define a
\emph{polyhedral cone} to be a finite union of \cpc s.

We recall two classical lemmas:

\begin{lemma}\label{L:Declpol}
Let $n$ be a positive integer, let $f\in\FA(n)$, let $K$ be a
polyhedral cone of $\QQ^n$. Then there exist a positive integer $N$ and a
decomposition $K=\bigcup_{1\leq l\leq N}K_l$ of $K$ into \cpc s such that
for all $l$, there exists a linear functional $f_l$ on $\ZZ^n$ such that
$f\res_{K_l}=f_l\res_{K_l}$.
\end{lemma}

\begin{proof}
Without loss of generality, $K$ is a \cpc.
We decompose $f$ as in \eqref{Eq:Decomp}. For every permutation $\sigma$
of $\fs(k+l)$ that leaves both sets $\{1,\ldots,k\}$ and
$\{k+1,\ldots,k+l\}$ invariant, we define $K_\sigma$ as the set of
those $x\in K$ such that for all $i\in\fs(k+l)$ distinct from $k$ and
from $k+l$, the inequality $f_{\sigma(i)}(x)\leq f_{\sigma(i+1)}(x)$ holds.
On each $K_\sigma$, $f$ is equal to some $f_i-f_j$. The conclusion of the
lemma follows, with $N=k!l!$.
\end{proof}

\begin{lemma}\label{L:pinvcpc}
Let $n$ be a positive integer, let $f\in\FA(n)$. Then the set
 \[
 \{x\in\QQ^n\mid f(x)\geq0\}
 \]
is a polyhedral cone of $\QQ^n$.\qed
\end{lemma}

If $K$ is a polyhedral cone of $\QQ^n$, a map $f\colon K\to\QQ$ is
\emph{piecewise linear}, if there exists a decomposition
$K=\bigcup_{1\leq l\leq N}K_l$ of $K$ into \cpc s such that
for all $l$, there exists a linear functional $f_l$ on $\QQ^n$ which
satisfies $f\res_{K_l}=f_l\res_{K_l}$. By Lemma~\ref{L:Declpol}, every
element of $\FA(n)$ is piecewise linear on $\QQ^n$. We denote by $\PL(K)$
the  \alo\ group of all piecewise linear maps from a polyhedral cone $K$
to $\QQ$.

The next lemma is folklore. It can also be easily deduced from
Proposition~\ref{P:GEHSeqsys2}.

\begin{lemma}\label{L:cpcPosGen}
Let $n$ be a positive integer. Then every \cpc\ $K$ of $\QQ^n$ can be
positively generated by a finite subset $X$, that is,
 \begin{equation}
 K=\left\{\sum_{x\in X}\lambda_xx\mid
 \lambda_x\in\QQ^+,\text{ for all }x\in X\right\}.\tag*{\qed}
 \end{equation}
\end{lemma}

\subsection{Finitely presented \alo\ groups}

An \alo\ group $\bG$ is \emph{finitely presented}, if it is isomorphic to
a quotient $\FA(n)/I$, where $n$ is a positive integer and
$I$ is a finitely generated $\ell$-ideal of
$\FA(n)$ (an \emph{$\ell$-ideal} is a convex subgroup closed under
$\wedge$ and $\vee$). If $I$ is generated by $p_1$,\dots, $p_m$, then $I$
is also generated by the single element $p=\sum_{i=1}^m|p_i|$ (where
$|x|=x\vee(-x)$ for all $x$). Hence, $I$ takes on the following simple
form
 \[
 I=\{f\in\FA(n)\mid\exists k\in\NN,\ |f|\leq kp\}.
 \]
Put $Z(f)=\{x\in\QQ^n\mid f(x)=0\}$, for all $f\in\FA(n)$. By
Lemma~\ref{L:pinvcpc}, applied to $-|f|$, $Z(f)$ is a polyhedral cone of
$\QQ^n$. Put $K=Z(p)$.

We define a sub-lattice-ordered group $\bH$ of $\PL(K)$ by
 \[
 \bH=\{f\res_K\mid f\in\FA(n)\}.
 \]
There is a natural surjective homomorphism of lattice-ordered groups,
$\rho\colon\bG\twoheadrightarrow\bH$,
defined by the rule $\rho(f+I)=f\res_K$, for all
$f\in\FA(n)$. By \cite[Lemma~3.3]{Bake}, $\rho$ is an \emph{isomorphism}.
We have thus outlined part of the proof of the following classical result:

\begin{lemma}\label{L:FPlgroup}
Every finitely presented \alo\ group can be embedded into the
lattice-ordered group $\PL(K)$ of all piecewise linear maps on a
polyhedral cone $K$.\qed
\end{lemma}

\begin{proof}[Proof of Theorem~\textup{\ref{T:FinPres}}]
By Lemma~\ref{L:FPlgroup}, in order to prove that every finitely presented
\alo\ group is \po-coherent, it suffices to prove that all the groups
$\PL(K)$ are \po-coherent. Thus let $n$ be a positive integer, let $K$ be a
polyhedral cone of $\QQ^n$, let $f_1$,\dots, $f_m\in\PL(K)$. We consider
the set
 \[
 \SS=\left\{\vv<\lambda_1,\ldots,\lambda_m>\in(\QQ^+)^m\mid
 \forall x\in K,\ \sum_{i=1}^m\lambda_if_i(x)\geq0\right\},
 \]
we must prove that $\SS\cap\ZZ^m$ is a finitely generated submonoid of
$\ZZ^m$. Since the $f_i$ are piecewise linear, there exist $N\in\NN$ and a
decomposition $K=\bigcup_{1\leq l\leq N}K_l$ of $K$ into \cpc s
such that every $f_i$ is
linear on every $K_l$, that is, there exists a linear functional $f_{i,l}$
on $\QQ^n$ such that $f_i\res_{K_l}=f_{i,l}\res_{K_l}$. Hence, the equality
 \[
 \SS=\bigcap_{1\leq l\leq N}\SS_l
 \]
holds, where we put
 \[
 \SS_l=\left\{\vv<\lambda_1,\ldots,\lambda_m>\in(\QQ^+)^m\mid
 \forall x\in K_l,\ \sum_{i=1}^m\lambda_if_{i,l}(x)\geq0\right\},
 \]
for all $l\in\fs(N)$. However, $K_l$ is positively generated by a finite
subset $X_l$ (see Lemma~\ref{L:cpcPosGen}), so that the equality
 \[
 \SS_l=\bigcap_{x\in X_l}
 \left\{\vv<\lambda_1,\ldots,\lambda_m>\in(\QQ^+)^m\mid
 \sum_{i=1}^m\lambda_if_{i,l}(x)\geq0\right\}
 \]
holds. In particular, $\SS_l$ is a \cpc, for all $l\in\fs(N)$, thus $\SS$
is a \cpc. This implies that the elements of $\SS\cap\ZZ^n$ are the
solutions in $\ZZ^n$ of a finite system of inequalities with coefficients
from $\QQ$, thus, after having cleared away denominators, from $\ZZ$. By
Corollary~\ref{C:GEHScoh}, we obtain that $\SS\cap\ZZ^n$ is a finitely
generated monoid, which completes the proof.
\end{proof}

\begin{example}\label{Ex:lgroups}
This example shows that the assumption that $\bF$ is finitely
generated convex is essential in the statement of
Proposition~\ref{P:ConvQuot}.

For this, we consider the abelian lattice-ordered groups $\bG$ and
$\bH$ defined as follows. First, $\bG=\FA(2)$ is the free abelian
lattice-ordered group on two generators, $x$ and $y$, and
$\bH=\ZZ+\alpha\ZZ$, viewed as an ordered additive subgroup of $\RR$,
where $\alpha$ is any positive irrational number. Furthermore, let
$f\colon\bG\twoheadrightarrow\bH$ be the unique homomorphism of
lattice-ordered groups sending $x$ to $1$, and $y$ to $\alpha$. By
Theorem~\ref{T:FinPres}, $\bG$ is a \po-coherent \poag.

Furthermore, observe that $\ker f$ is an ideal (that is, a convex,
directed subgroup) of $\bG$, and that $\bH$ is isomorphic to
$\bG/\!\ker f$. However, $\bH$ is not a \po-coherent \porm{\ZZ}.
Indeed, otherwise, the following set
 \[
 \SS=\{\vv<x,y>\in(\ZZ^+)^2\mid x-\alpha y\geq0\}
 \]
would be a finitely generated submonoid of $(\ZZ^+)^2$. Let
 \[
 \{\vv<x_i,y_i>\mid 1\leq i\leq n\}
 \]
be a finite generating subset of
$\SS$, with $n\in\NN$ and $x_i>0$, for all $i$. Let $\lambda$ be the
largest element of $\{y_i/x_i\mid 1\leq i\leq n\}$. So
$\alpha\lambda\leq1$, thus, since $\alpha$ is irrational and $\lambda$ is
rational nonzero, $\alpha\lambda<1$. Hence, there exists $\vv<x,y>$ in
$\NN\times\NN$ such that $\alpha y/x<1$ and $\lambda<y/x$. The first
inequality implies that $\vv<x,y>\in\SS$, so there are nonnegative
integers $k_i$, for $1\leq i\leq n$, such that $x=\sum_ik_ix_i$ and
$y=\sum_ik_iy_i$. In particular, $y\leq\lambda x$, so $y/x\leq\lambda$, a
contradiction.
\end{example}

\section{Discussion}

From the results of all the previous sections, we can deduce the
following:

\begin{theorem}\label{T:SubMod}
Let $\bR$ be either a subring of $\QQ$ or a totally ordered
division ring. Let $\bA$ be a \porm{\bR}.
Then the following properties hold:

\begin{enumerate}
\item $\bA$ is finitely \po-presented if and only if it is finitely
related, if and only if $\bA$ is a finitely presented $\bR$-module and
$\bA^+$ is a finitely generated $\bR^+$-semimodule.

\item Suppose that $\bA$ is finitely \po-presented.
Then every finitely generated submodule of $\bA$ is finitely
\po-presented.
\end{enumerate}
\end{theorem}

\begin{proof}
By Theorem~\ref{T:CohSubQ} and Corollary~\ref{C:GEHScoh2}, $\bR$ is a
\po-coherent \por. We conclude the proof by Theorem~\ref{T:FinPresEquiv}
for (i), and by Theorem~\ref{T:SubmFP} for (ii).
\end{proof}

Since $\ZZ$ is a n\oe therian ring, the particular case where $\bR=\ZZ$
states that (i) a \poag\ $\bG$ is finitely \po-presented if and only if
$\bG$ is a finitely generated group and $\bG^+$ is a finitely generated
monoid, and (ii) every subgroup of a finitely \po-presented \poag\ is
finitely \po-presented. These results are established by completely
different means in \cite{CaWe}.

\begin{problem}\label{Pb:fpfr}
Let $\bR$ be a \por, let $\bA$ be a \porm{\bR}. Suppose that $\bA$ is
finitely \po-presented. Is $\bA$ finitely related?
\end{problem}

The converse, namely the fact that finitely related implies
finitely \po-presented, holds as a rule, see Proposition~\ref{P:FRimplFP}.
Furthermore, for partially ordered right modules over right \po-coherent
rings, finitely \po-presented is equivalent to finitely related, see
Theorem~\ref{T:FinPresEquiv}. On the other hand, it seems unlikely that
the answer to Problem~\ref{Pb:fpfr} could be positive in general.

\begin{problem}\label{Pb:CohReal}
Let $\bR$ be a subring of the real field $\RR$. Suppose that
$\bR$ is coherent. Is $\bR$ \po-coherent?
More generally, if $\bR$ is a totally ordered ring,
is it the case that coherence of $\bR$ implies \po-coherence of $\bR$?
\end{problem}

We do not know the answer even in the case where $\bR$ is a n\oe therian
subring of $\RR$, even in the case where $\bR$ is an extension
of $\ZZ$ by finitely many real algebraic numbers. Concretely, the
problem is the following: given a coherent subring $\bR$ of $\RR$, and
$n\in\NN$ and $\alpha_1$, \dots, $\alpha_n\in\bR$, prove that the
solution set in $(\bR^+)^n$ of the equation
 \[
 \alpha_1x_1+\cdots+\alpha_nx_n=0
 \]
is a finitely generated $\bR^+$-subsemimodule of $(\bR^+)^n$.

\begin{problem}\label{Pb:CohReal2}
Let $\bR$ be a subring of the real field, totally ordered with its
natural ordering. Suppose that $\bR$ is \po-coherent. Does there exist an
analogue of the Grillet, Effros, Handelman, and Shen Theorem that would
hold for \porm{\bR}s?
\end{problem}

\begin{problem}
Investigate right \po-coherence of subrings of totally ordered division
rings.
\end{problem}

\section*{Acknowledgment}
The author would like to thank the anonymous referee for having thoroughly
read the paper, pointing an embarrassing number of errors and oversights,
and suggesting a number of improvements.

\end{document}